\documentclass[smallextended]{svjour3}       % onecolumn (second format)
\smartqed  % flush right qed marks, e.g. at end of proof
\usepackage{graphicx}
\usepackage{latexsym}
\usepackage{calc}
\usepackage{amssymb, amsmath, amsfonts, listings, mathrsfs}
\usepackage{hyperref}

\spnewtheorem{assumption}{Assumption}{\bf}{\it}

\newcommand{\Def}{\stackrel{\mathrm{def}}{=}}

\newcommand{\dom}{{\rm dom \,}}
\newcommand{\beq}{\begin{equation}}
\newcommand{\eeq}{\end{equation}}

\newcommand{\R}{\mathbb{R}}

\newcommand{\E}{\mathbb{E}}
\newcommand{\Q}{\mathcal{Q}}
\newcommand{\G}{\mathcal{G}}

\newcommand{\SetEQ}{\setcounter{equation}{0}}
\newcommand{\refLE}[1]{\ensuremath{\stackrel{(\ref{#1})}{\leq}}}
\newcommand{\refEQ}[1]{\ensuremath{\stackrel{(\ref{#1})}{=}}}
\newcommand{\refGE}[1]{\ensuremath{\stackrel{(\ref{#1})}{\geq}}}

\newcommand{\ba}{\begin{array}}
\newcommand{\ea}{\end{array}}
\newcommand{\beann}{\begin{eqnarray*}}
\newcommand{\eeann}{\end{eqnarray*}}
\newcommand{\bea}{\begin{eqnarray}}
\newcommand{\eea}{\end{eqnarray}}

\newcommand{\BT}{\begin{theorem}}
\newcommand{\ET}{\end{theorem}}
\newcommand{\BL}{\begin{lemma}}
\newcommand{\EL}{\end{lemma}}
\newcommand{\BC}{\begin{corollary}}
\newcommand{\EC}{\end{corollary}}
\newcommand{\BE}{\begin{example}}
\newcommand{\EE}{\end{example}}
\newcommand{\BD}{\begin{definition}}
\newcommand{\ED}{\end{definition}}
\newcommand{\BR}{\begin{remark}}
\newcommand{\ER}{\end{remark}}
\newcommand{\BAS}{\begin{assumption}}
\newcommand{\EAS}{\end{assumption}}
\newcommand{\BI}{\begin{itemize}}
\newcommand{\EI}{\end{itemize}}

\newcommand{\BMP}{\begin{minipage}{9.5cm}}
\newcommand{\EMP}{\end{minipage}}
\newcommand{\MPT}{\begin{minipage}{11.5cm}}
\newcommand{\EPT}{\end{minipage}}

\newcommand{\la}{\langle}
\newcommand{\ra}{\rangle}

\begin{document}

\title{Local convergence of
	   tensor methods\thanks{The research results of this paper were obtained
	   	in the framework of ERC Advanced Grant 788368.}
}


\author{Nikita Doikov         \and
        Yurii Nesterov
}


\institute{Nikita Doikov \at
              Institute of Information and Communication Technologies, Electronics and Applied Mathematics (ICTEAM), Catholic
              University of Louvain (UCL),  Louvain-la-Neuve, Belgium. \\
               \email{Nikita.Doikov@uclouvain.be}. 
               ORCID: 0000-0003-1141-1625.
           \and
           Yurii Nesterov \at
              Center for Operations Research and Econometrics (CORE), Catholic University of Louvain
              (UCL), 34 voie du Roman Pays, 1348 Louvain-la-Neuve, Belgium. \\
              \email{Yurii.Nesterov@uclouvain.be}. 
              ORCID: 0000-0002-0542-8757.
}

\date{Received: 16 December 2019 / Accepted: 8 December 2020 \\
	\textcopyright\,\!\!  The Author(s) 2021}

\maketitle

\begin{abstract}
In this paper, we study local convergence of
high-order Tensor Methods for solving convex optimization
problems with composite objective. We justify local
superlinear convergence under the assumption of uniform
convexity of the smooth component, having
Lipschitz-continuous high-order derivative. The
convergence both in function value and in the norm of
minimal subgradient is established. Global complexity
bounds for the Composite Tensor Method in convex and
uniformly convex cases are also discussed. Lastly, we show
how local convergence of the methods can be globalized
using the inexact proximal iterations.
\end{abstract}

\section{Introduction}
\SetEQ

\vspace{1ex}\noindent
{\bf Motivation.} In Nonlinear Optimization, it seems to
be a natural idea to increase the performance of numerical
methods by employing high-order oracles. However, the main
obstacle to this approach consists in a prohibiting
complexity of the corresponding Taylor approximations
formed by the high-order multidimensional polynomials,
which are difficult to store, handle, and minimize. If we
go just one step above the commonly used quadratic
approximation, we get a multidimensional polynomial of
degree three which is never convex. Consequently, its
usefulness for optimization methods is questionable.

However, recently in \cite{nesterov2019implementable}
it was shown that the Taylor polynomials of
{\em convex functions} have a very interesting structure.
It appears that their augmentation by a power of Euclidean
norm with a reasonably big coefficients gives us a global
upper {\em convex} model of the objective function,
which keeps all advantages of the local high-order
approximation.

One of the classical and well-known results in Nonlinear
Optimization is related to the local quadratic convergence
of Newton's
method~\cite{kantorovich1948functional,nesterov2018lectures}.
Later on, it was generalized to the  case of
\textit{composite} optimization
problems~\cite{lee2014proximal}, where the objective is
represented as a sum of two convex components: smooth, and
possibly nonsmooth but simple. Local superlinear convergence
of the Incremental Newton method for
finite-sum minimization problems
was established in~\cite{rodomanov2016superlinearly}.

The study of high-order numerical methods for solving
nonlinear equations is dated back to the work
of Chebyshev in 1838,
where the scalar methods of order three 
and four were proposed~\cite{chebyshev1951sobranie}.
The methods of arbitrary order for solving
nonlinear equations were studied in~\cite{evtushenko2014methods}.

A big step in the second-order optimization theory was
made since~\cite{nesterov2006cubic}, where Cubic
regularization of the Newton method with its global
complexity estimates was proposed. Additionally, the local
superlinear convergence was justified.
See also~\cite{cartis2011adaptive1} for the local analysis
of the Adaptive cubic regularization methods.

Our paper is aimed to study local convergence of
high-order methods, generalizing corresponding results
from~\cite{nesterov2006cubic} in several ways. We
establish local superlinear convergence of Tensor
Method~\cite{nesterov2019implementable} of degree $p \geq
2$, in the case when the objective is composite, and its
smooth part is uniformly convex of arbitrary degree $q$
from the interval $2 \leq q < p + 1$. For strongly convex
functions ($q=2$), this gives the local convergence of
degree $p$.

\vspace{1ex}\noindent
{\bf Contents.}
We formulate our problem of interest and define
a step of the Regularized Composite Tensor Method
in Sect.~\ref{sc-MIn}.
Then, we declare some of its properties, which are
required for our analysis.

In Sect.~\ref{sc-LocF}, we prove local superlinear
convergence of the Tensor Method in function value, and in
the norm of minimal subgradient, under the assumption of
uniform convexity of the objective.

In Sect.~\ref{sc-GlobG}, we discuss global behavior of
the method and justify sublinear and linear global rates
of convergence for convex and uniformly convex cases,
respectively.

One application of our developments is provided in
Sect.~\ref{sc-Prox}. We show how local convergence can
be applied for computing an inexact step in proximal
methods. A global sublinear rate of convergence for the
resulting scheme is also given.

\vspace{1ex}\noindent
{\bf Notations and generalities.} In what follows, we
denote by $\E$ a finite-dimensional real vector space, and
by $\E^*$ its dual spaced composed by linear functions on
$\E$. For such a function $s \in \E^*$, we denote by $\la
s, x \ra$ its value at $x \in \E$. Using a self-adjoint
positive-definite operator $B: \E \to \E^*$ (notation $B =
B^* \succ 0$), we can endow these spaces with mutually
conjugate Euclidean norms:
$$
\ba{rcl}
\| x \| & = & \la B x, x \ra^{1/2}, \quad x \in \E, \quad
\| g \|_* \; = \; \la g, B^{-1} g \ra^{1/2}, \quad g \in
\E^*.
\ea
$$

For a smooth function $f: \dom f \to \R$ with convex and
open domain $\dom f  \subseteq \E$, denote by $\nabla
f(x)$ its gradient, and by $\nabla^2 f(x)$ its Hessian
evaluated at point $x \in \dom f \subseteq \E$. Note that
$$
\ba{rcl}
\nabla f(x) & \in & \E^*, \quad \nabla^2 f(x) h \; \in \;
\E^*, \quad x \in \dom f, \; h \in \E.
\ea
$$
For non-differentiable convex function $f(\cdot)$, we
denote by $\partial f(x) \subset \E^*$ its subdifferential
at the point $x \in \dom f$.

In what follows, we often work with directional
derivatives. For $p \geq 1$, denote by
$$
\ba{c}
D^p f(x)[h_1, \dots, h_p]
\ea
$$
the directional derivative of function $f$ at $x$ along
directions $h_i \in \E$, $i = 1, \dots, p$. If all
directions $h_1, \dots, h_p$ are the same, we apply a
simpler notation
$$
\ba{c}
D^p f(x)[h]^p, \quad h \in \E.
\ea
$$

Note that $D^p f(x)[ \cdot]$ is a {\em symmetric
	$p$-linear form}. Its {\em norm} is defined in the
standard way:
\beq\label{eq-DNorm1}
\ba{rcl}
\| D^pf(x) \| & = & \max\limits_{h_1, \dots, h_p \in \E}
\left\{ D^p f(x)[h_1, \dots, h_p ]: \; \| h_i \| \leq 1,
\, i = 1,
\dots, p \right\}\\
\\
& = & \max\limits_{h \in \E} \left\{ \Big| D^p
f(x)[h]^p\Big|: \; \| h \| \leq 1 \right\}
\ea
\eeq
(for the last equation see, for example, Appendix 1 in
\cite{nesterov1994interior}). Similarly, we define
\beq\label{eq-DNorm2}
\ba{rcl}
\| D^pf(x) - D^pf(y) \| & = & \max\limits_{h \in \E}
\left\{ \Big| D^p f(x)[h]^p - D^pf(y)[h]^p\Big|: \; \| h
\| \leq 1 \right\}.
\ea
\eeq

In particular, for any $x \in \dom f$ and $h_1, h_2 \in
\E$, we have
$$
\ba{rcl}
Df(x)[h_1] & = & \la \nabla f(x), h_1 \ra, \quad
D^2f(x)[h_1, h_2] \; = \; \la \nabla^2 f(x) h_1, h_2 \ra.
\ea
$$
Thus, for the Hessian, our definition corresponds to a
{\em spectral norm} of the self-adjoint linear operator
(maximal module of all eigenvalues computed with respect
to $B \succ 0$).

Finally, the Taylor approximation of function $f(\cdot)$
at $x \in \dom f$ is defined as follows:
$$
\ba{rcl}
f(x+h) & = & \Omega_p(f, x; x + h)  + o(\|h\|^p), \quad x+h \in
\dom
f,\\
\\
\Omega_p(f,x;y) & \Def & f(x) + \sum\limits_{k=1}^p {1
	\over k!} D^k f(x)[y-x]^k, \quad y \in \E.
\ea
$$
Consequently, for all $y \in \E$ we have
\beq\label{eq-PGrad}
\ba{rcl}
\nabla \Omega_p(f,x;y) & = & \sum\limits_{k=1}^p {1 \over
	(k-1)!} D^k f(x)[y-x]^{k-1},
\ea
\eeq
\beq\label{eq-PHess}
\ba{rcl}
\nabla^2 \Omega_p(f,x;y) & = & \sum\limits_{k=2}^p {1
	\over (k-2)!} D^k f(x)[y-x]^{k-2}.
\ea
\eeq

\section{Main inequalities}\label{sc-MIn}
\SetEQ

In this paper, we consider the following {\em composite}
convex minimization problem
\beq\label{prob-Main}
\min\limits_{x \in \dom h} \Big\{ F(x) = f(x) + h(x)
\Big\},
\eeq
where $h: \E \to \R \cup \{+\infty\}$ is a {\em simple} proper
closed convex function and $f \in C^{p,p}(\dom h)$ for a
certain $p \geq 2$. In other words, we assume that the
$p$th derivative of function $f$ is Lipschitz continuous:
\beq\label{eq-Lip}
\ba{rcl}
\| D^p f(x) - D^p f(y) \| & \leq & L_p \| x - y \|,
\quad x, y \in \dom h.
\ea
\eeq

Assuming that $L_{p} < +\infty$, by the standard
integration arguments we can bound the residual between
function value and its Taylor approximation:
\beq\label{eq-BoundF}
\ba{rcl}
| f(y) - \Omega_p(f,x;y) | & \leq & {L_{p} \over (p+1)!}
\| y - x \|^{p+1}, \quad x, y \in \dom h.
\ea
\eeq
Applying the same reasoning to functions $\la \nabla
f(\cdot), h \ra$ and $\la \nabla^2 f(\cdot) h, h \ra$ with
direction $h \in \E$ being fixed, we get the following
guarantees:
\beq\label{eq-BoundG}
\ba{rcl}
\| \nabla f(y) - \nabla \Omega_p(f,x;y) \|_* & \leq & {L_p
	\over p!} \| y - x \|^{p},
\ea
\eeq
\beq\label{eq-BoundH}
\ba{rcl}
\| \nabla^2 f(y) - \nabla^2 \Omega_p(f,x;y) \| & \leq &
{L_p \over (p-1)!} \| y - x \|^{p-1},
\ea
\eeq
which are valid for all $x, y \in \dom h$.

Let us define now one step of the {\em Regularized
	Composite Tensor Method} (RCTM) of degree $p \geq 2$:
\beq\label{def-RCTS}
\ba{rcl}
T & \equiv & T_H(x) \; \Def \; \arg\min\limits_{y \in \E}
\left\{ \Omega_p(f,x;y) + {H \over (p+1)!} \| y - x
\|^{p+1} + h(y) \right\}.
\ea
\eeq
It can be shown that for
\beq\label{eq-H}
\ba{c}
\mbox{\fbox{\rule[-2mm]{0cm}{6mm}$H \; \geq \; p L_p$}}
\ea
\eeq
the auxiliary optimization problem in (\ref{def-RCTS}) is
{\em convex} (see Theorem 1 in \cite{nesterov2019implementable}).
This condition
is crucial for implementability of our methods and we
always assume it to be satisfied.

Let us write down the first-order optimality condition for
the auxiliary optimization problem in (\ref{def-RCTS}):
\beq\label{eq-OptC}
\ba{rcl}
\la \nabla \Omega_p(f,x;T) + {H \over p!} \| T - x
\|^{p-1}B(T-x), y - T \ra + h(y) & \geq & h(T),
\ea
\eeq
for all $y \in \dom h$.
In other words, for vector
\beq\label{def-HP}
\ba{rcl}
h'(T) & \Def & - \left( \nabla \Omega_p(f,x;T) + {H \over
	p!} \| T - x \|^{p-1}B(T-x) \right)
\ea
\eeq
we have $h'(T) \stackrel{(\ref{eq-OptC})}{\in} \partial
h(T)$. This fact explains our notation
\beq\label{def-FP}
\ba{rcl}
F'(T) & \Def & \nabla f(T) + h'(T) \; \in \partial F(T).
\ea
\eeq

Let us present some properties of the point $T = T_H(x)$.
First of all, we need some bounds for the norm of vector
$F'(T)$. Note that
\beq\label{eq-N3}
\ba{rcl}
\Big\| F'(T) + {H \over p!} \| T - x \|^{p-1}B(T-x)
\Big\|_* & \refEQ{def-HP} & \Big\| \nabla f(T) - \nabla
\Omega_p(f,x;T) \Big\|_*\\
\\
& \refLE{eq-BoundG} & {L_p \over p!} \| T - x \|^p.
\ea
\eeq
Consequently,
\beq\label{eq-NewG}
\ba{rcl}
\| F'(T) \|_* & \leq & {L_p+H \over p!} \| T - x \|^p.
\ea
\eeq
Secondly, we use the following lemma. 

\BL\label{lm-DecF}
Let $\beta > 1$ and $H = \beta L_p$. Then
\beq\label{eq-DecFB}
\ba{rcl}
\la F'(T), x - T \ra & \geq & \left( {p! \over (p+1)L_p}
\right)^{1 \over p} \cdot \| F'(T) \|_*^{p+1 \over p}
\cdot {(\beta^2 - 1)^{p-1 \over 2p} \over \beta} \cdot {p
	\over (p^2-1)^{p-1 \over 2p}}.
\ea
\eeq
In particular, if $\beta = p$, then
\beq\label{eq-DecF}
\ba{rcl}
\la F'(T), x - T \ra & \geq & \left( {p! \over (p+1)L_p}
\right)^{1 \over p} \cdot \| F'(T) \|_*^{p+1 \over p}.
\ea
\eeq
\EL
\proof
Denote $r = \| T  - x \|$, $h = {H \over p!}$, and $l =
{L_p \over p!}$. Then inequality (\ref{eq-N3}) can be
written as follows:
$$
\ba{rcl}
\| F'(T) + h r^{p-1} B(T-x) \|^2_* & \leq & l^2 r^{2p}.
\ea
$$
This means that
\beq\label{eq-N4}
\ba{rcl}
\la F'(T), x - T \ra & \geq & {1 \over 2 h r^{p-1}} \|
F'(T) \|_*^2 + {r^{2p} (h^2 - l^2) \over 2h r^{p-1}}.
\ea
\eeq
Denote
$$
\ba{rcl}
a & = & {1 \over 2h} \| F'(T) \|_*^2, \quad b \; = \; {h^2
	- l^2 \over 2h}, \quad \tau \; = \; r^{p-1}, \quad \alpha
\; = \; {p+1 \over p-1}.
\ea
$$
Then inequality (\ref{eq-N4}) can be rewritten as follows:
$$
\ba{rcl}
\la F'(T) , x - T \ra & \geq & {a \over \tau} + b
\tau^{\alpha} \; \geq \; \min\limits_{t > 0} \left\{ {a
	\over t} + b t^{\alpha} \right\} \; = \; (1+\alpha)
\left({a \over \alpha} \right)^{\alpha \over 1 + \alpha}
b^{1 \over 1 + \alpha}.
\ea
$$
Taking into account that $1+\alpha = {2p \over p-1}$ and
${\alpha \over 1 + \alpha} = {p + 1 \over 2p}$, and using the
actual meaning of $a$, $b$, and $\alpha$, we get
$$
\ba{rcl}
\la F'(T), x - T \ra & \geq & {2 p \over p-1} \cdot { \|
	F'(T) \|_*^{p+1 \over p} \over (2h)^{p+1 \over 2p}} \cdot
{(p-1)^{p+1 \over 2p} \over (p+1)^{p+1 \over 2p}} \cdot
{(h^2 - l^2)^{p-1 \over 2p} \over (2h)^{p-1 \over 2p}}\\
\\
& = & \| F'(T) \|_*^{p+1 \over p} \cdot {(h^2 - l^2)^{p-1
		\over 2p} \over h} \cdot {p \over (p+1)^{p+1 \over 2p}
	(p-1)^{p-1 \over 2p}}\\
\\
& = & \| F'(T) \|_*^{p+1 \over p} \cdot {(h^2 - l^2)^{p-1
		\over 2p} \over h} \cdot {p \over (p^2-1)^{p-1 \over 2p}
	(p+1)^{1 \over p}}.
\ea
$$
It remains to note that
$$
\ba{rcl}
{(h^2 - l^2)^{p-1 \over 2p} \over h} & = & {(H^2 -
	L_p^2)^{p-1 \over 2p} \over H} \cdot (p!)^{1 \over p} \; =
\; {(\beta^2 - 1)^{p-1 \over 2p} \over \beta} \cdot
\left({p! \over L_p} \right)^{1 \over p}.
\ea
$$
\qed

\section{Local convergence}\label{sc-LocF}
\SetEQ

The main goal of this paper consists in analyzing the
local behavior of the {\em Regularized Composite Tensor
	Method} (RCTM):
\beq\label{met-RCTM}
\ba{rcl}
x_0 \; \in \; \dom h, \quad x_{k+1} & = & T_H(x_k), \quad
k \geq 0,
\ea
\eeq
as applied to the problem (\ref{prob-Main}). In order to
prove local superlinear convergence of this scheme, we
need one more assumption.

\begin{assumption}\label{assumption-Uni}
	The objective in problem
	(\ref{prob-Main}) is uniformly convex of degree $q \geq 2$.
	Thus, for all $x, y \in \dom h$ and for all $G_x \in \partial F(x), G_y \in \partial F(y)$,
	it holds:
	\beq\label{eq-UniC}
	\ba{rcl}
	\la G_x - G_y, x - y \ra & \geq & \sigma_q
	\| x - y \|^q, 
	\ea
	\eeq
	for certain $\sigma_q > 0$.
\end{assumption}

It is well known that this assumption guarantees the
uniform convexity of the objective function
(see, for example, Lemma 4.2.1 in
\cite{nesterov2018lectures}):
\beq\label{eq-UniF}
\ba{rcl}
F(y) & \geq & F(x) + \la G_x, y - x \ra + {\sigma_q \over
	q} \| y - x \|^q, \quad y \in \dom h,
\ea
\eeq
where $G_x$ is an arbitrary subgradient from $\partial
F(x)$. Therefore,
\beq\label{eq-DecFG}
\ba{rcl}
F^* & = & \min\limits_{y \in \dom h} F(y) \; \geq \;
\min\limits_{y \in \E} \left\{ F(x) + \la G_x, y
- x \ra + {\sigma_q \over q} \| y - x \|^q \right\}\\
\\
& = & F(x) - {q-1 \over q} \left({1 \over \sigma_q}
\right)^{1 \over q-1} \| G_x \|_*^{q \over q-1}.
\ea
\eeq
This simple inequality gives us the following local
convergence rate for RCTM.
\BT\label{th-LocF}
For any $k \geq 0$ we have
\beq\label{eq-RateF}
\ba{rcl}
F(x_{k+1}) - F^* & \leq & (q-1) q^{p-q+1 \over q-1} 
\bigl({1 \over \sigma_q}\bigr)^{p+1 \over q-1} \left({L_p + H
	\over p!} \right)^{q \over q - 1} \big[F(x_k) - F^* \big]^{p
	\over q-1}.
\ea
\eeq
\ET
\proof
Indeed, for any $k \geq 0$ we have
$$
\ba{rcl}
F(x_k) - F^* & \geq & F(x_k) - F(x_{k+1}) \\
\\
& \refGE{eq-UniF} &  \la
F'(x_{k+1}), x_k - x_{k+1} \ra + {\sigma_q \over q} \| x_k - x_{k+1} \|^q\\
\\
& \refGE{eq-DecFB} & {\sigma_q \over q} \| x_k - x_{k+1}
\|^q \; \refGE{eq-NewG} \; {\sigma_q \over q} \left( {p!
	\over L_p+H} \| F'(x_{k+1}) \|_* \right)^{q
	\over p}\\
\\
& \refGE{eq-DecFG} & {\sigma_q \over q} \left( {p! \over
	L_p+H}\right)^{q \over p} \left( {q \, \sigma_q^{1 \over
		q-1} \over q-1} (F(x_{k+1})-F^*)
\right)^{q -1\over p}.
\ea
$$
And this is exactly inequality (\ref{eq-RateF}).
\qed

\BC\label{cor-LocF}
If $p > q-1$, then method (\ref{met-RCTM}) has local
superlinear rate of convergence for problem
(\ref{prob-Main}).
\EC
\proof
Indeed, in this case ${p \over q-1} > 1$.
\qed

For example, if $q = 2$ (strongly convex function) and
$p=2$ (Cubic Regularization of the Newton Method), then
the rate of convergence is quadratic. If $q=2$, and $p =
3$, then the local rate of convergence is cubic, etc.

Let us study now the local convergence of the method
(\ref{met-RCTM}) in terms of the norm of gradient. For any
$x \in \dom h$ denote
\beq\label{def-GNorm}
\ba{rcl}
\eta(x) & \Def & \min\limits_{g \in \partial h(x)} \| \nabla
f(x) + g \|_*.
\ea
\eeq
If $\partial h(x) = \emptyset$, we set $\eta(x) =
+\infty$.

\BT\label{th-RateG}
For any $k \geq 0$ we have
\beq\label{eq-RateG}
\ba{rcl}
\eta(x_{k+1})
& \, \leq \, &
\|F'(x_{k + 1}) \|_{*}
\; \leq \;
{L_p + H \over p!} \left[ {1 \over
	\sigma_q} \, \eta(x_k) \right]^{p \over q-1}.
\ea
\eeq
\ET
\proof
Indeed, in view of inequality \eqref{eq-UniC}, we have
$$
\ba{rcl}
\la \nabla f(x_k) + g_k, x_{k} - x_{k + 1} \ra
& \geq &
\la F'(x_{k + 1}), x_k - x_{k + 1} \ra + \sigma_q \|x_k - x_{k + 1}\|^q \\
\\
& \refGE{eq-DecFB} & \sigma_q \|x_k - x_{k + 1}\|^{q},
\ea
$$
where $g_k$ is an arbitrary vector from $\partial h(x_k)$.
Therefore, we conclude that
$$
\ba{rcl}
\eta(x_k)  & \geq &  \sigma_q \| x_k -
x_{k+1} \|^{q-1}.
\ea
$$
It remains to use inequality
(\ref{eq-NewG}).
\qed

As we can see, the condition for superlinear convergence
of the method (\ref{met-RCTM}) in terms of the norm of the
gradient is the same as in Corollary \ref{cor-LocF}: we
need to have ${p \over q-1} > 1$, that is $p > q-1$.
Moreover, the local rate of convergence has the same order
as that for the residual of the function value.

According to Theorem~\ref{th-LocF}, the region of
superlinear convergence of RCTM in terms of the function
value is as follows:
\beq\label{eq-RegF}
\ba{rcl}
\Q & = & \left\{ x \in \dom h: \;  F(x) - F^* \; \leq \;
{1 \over q} \cdot \biggl(  { \sigma_q^{p + 1} \over (q -
	1)^{q - 1} } \cdot \Bigl( { p! \over L_p + H  } \Bigr)^{q}
\biggr)^{1 \over p - q + 1} \right\}.
\ea
\eeq
Alternatively, by Theorem~\ref{th-RateG}, in terms of the
norm of minimal subgradient~(\ref{def-GNorm}), the region
of superlinear convergence looks as follows:
\beq\label{eq-RegG}
\ba{rcl}
\G & = & \left\{ x \in \dom h: \; \eta(x) \; \leq \;
\biggl(  \sigma_q^{p} \cdot \Bigl( { p! \over L_p + H }
\Bigr)^{q - 1} \biggr)^{1 \over p - q + 1}  \right\}.
\ea
\eeq
Note that these sets can be very different. Indeed, set
$\Q$ is a closed and convex neighborhood of the point
$x^*$. At the same time, the structure of the set $\G$ can
be very complex since in general the function $\eta(x)$ is
discontinuous. Let us look at simple example where $h(x) =
\mbox{Ind}_Q(x)$, the indicator function of a closed
convex set $Q$.
\BE
Consider the following optimization problem:
\beq\label{prob-Eta}
\min\limits_{x \in \R^2} \left\{ f(x)  : \; \| x \|^2
\Def (x^{(1)})^2 + (x^{(2)})^2 \leq 1 \right\},
\eeq
with
$$
\ba{rcl}
f(x) & = &  
\frac{\sigma_2}{2}\|x - \bar{x}\|^2
+
\frac{2 \sigma_3}{3}\|x - \bar{x}\|^3,
\ea
$$
for some fixed $\sigma_2, \sigma_3 > 0$
and $\bar{x} = (0, -2) \in \R^2$.
We have
$$
\ba{rcl}
\nabla f(x) = r(x) \cdot ( x^{(1)}, x^{(2)} + 2),
\ea
$$
where $r: \R^2 \to \R$ is
$$
\ba{rcl}
r(x) = \sigma_2 + 2\sigma_3 \|x - \bar{x}\|.
\ea
$$
Note that $f$ is uniformly convex
of degree $q = 2$ with constant $\sigma_2$, 
and for $q = 3$ with constant $\sigma_3$
(see Lemma~4.2.3 in~\cite{nesterov2018lectures}).
Moreover, we have for any $\nu \in [0, 1]$:
$$
\ba{rcl}
\la \nabla f(x) - \nabla f(y), x - y \ra
& \geq & \sigma_2 \|x - y\|^2 + \sigma_3 \|x - y\|^3 \\
\\
& \geq & \min\limits_{t \geq 0} 
\Bigl\{ \frac{\sigma_2}{t^{\nu}} + \sigma_3 t^{1 - \nu}
\Bigr\} \cdot \|x - y\|^{2 + \nu} \\
\\
& \geq & \sigma_2^{1 - \nu} \sigma_3^{\nu} \cdot \|x - y\|^{2 + \nu}.
\ea
$$
Hence, this function is uniformly convex of any
degree $q \in [2, 3]$. At the same time,
the Hessian of $f$ is Lipschitz continuous
with constant $L_2 = 4 \sigma_3$
(see Lemma~4.2.4 in~\cite{nesterov2018lectures}).

Clearly, 
in this problem $x^*=(0,-1)$,
and it can be written in the composite form (\ref{prob-Main}) with
$$
\ba{rcl}
h(x) & = & \left\{ \ba{rl} +
\infty, & \mbox{if $\| x \| > 1$,} \\ 0,
&\mbox{otherwise.} \ea \right.
\ea
$$
Note that for $x \in \dom h \equiv \{ x: \; \| x \| \leq 1\}$, we have
$$
\ba{rcl}
\partial h(x) \; = \; \left\{ \ba{cl} 0, & \mbox{if $\|
	x \| < 1$,} \\ \{ \gamma x, \, \gamma \geq 0 \}, &\mbox{if
	$\| x \| = 1$.}
\ea \right.
\ea
$$
Therefore, if $\| x \| < 1$, then $\eta(x) = \| \nabla
f(x) \| \geq \sigma_2$. If $\| x \| = 1$, then
$$
\ba{rcl}
\eta^2(x) & \refEQ{def-GNorm} & 
\min\limits_{\gamma \geq
	0}
\Bigl\{ 
\bigl[ (r(x) + \gamma) x^{(1)} \bigr]^2
+ 
\bigl[ (r(x) + \gamma) x^{(2)} + 2 r(x) \bigr]^2
\Bigr\} \\
\\
& = &
\min\limits_{\gamma \geq
	0}
\Bigl\{ 
(r(x) + \gamma)^2 + 4r(x) (r(x) + \gamma) x^{(2)}  
+ 4 r^2(x) 
\Bigr\} \\
\\
& = &
\left\{ 
\ba{cl}   
	4r^2(x) (1 - (x^{(2)})^2), & \mbox{if $x^{(2)} \leq -\frac{1}{2}$,} \\ 
	r^2(x) (5 + 4 x^{(2)}), & \mbox{otherwise.}
\ea \right.
\ea
$$
Thus, in any neighbourhood
of $x^*$, $\eta(x)$ vanishes only along the boundary of
the feasible set.
\qed
\EE

So, the question arises how the Tensor Method
(\ref{met-RCTM}) could come to the region $\G$. The answer
follows from the inequalities derived in Section
\ref{sc-MIn}. Indeed,
$$
\ba{rcl}
\| F'(x_{k+1}) \|_* & \refLE{eq-NewG} & {L_p + H \over p!}
\| x_k - x_{k+1} \|^p,
\ea
$$
and
$$
\ba{rcl}
F(x_k) - F(x_{k+1}) & \geq & \la F'(x_{k+1}), x_k -
x_{k+1} \ra \\
\\
& \refGE{eq-DecF} &
\left( {p! \over (p+1)L_p} \right)^{1 \over p} \cdot \|
F'(x_{k+1}) \|^{p+1 \over p}_*.
\ea
$$
Thus, at some moment the norm $\| F'(x_k) \|_*$ will be
small enough to enter $\G$.

\section{Global complexity bounds}\label{sc-GlobG}
\SetEQ

Let us briefly discuss the global complexity bounds of the
method~(\ref{met-RCTM}), namely the number of iterations
required for coming from an arbitrary initial point $x_0
\in \dom h$ to the region~$\Q$. First, note that for every
step $T = T_H(x)$ of the method with parameter $H \geq p
L_p$, we have
$$
\ba{rcl}
F(T) & \refLE{eq-BoundF} &
\Omega_p(f,x;T) + \frac{H}{(p + 1)!}\|T - x\|^{p + 1} + h(T) \\
\\
& \refEQ{def-RCTS} & \min\limits_{y \in \E} \Bigl\{
\Omega_p(f,x;y) + \frac{H}{(p + 1)!}\|y - x\|^{p + 1} +
h(y)
\Bigr\} \\
\\
& \refLE{eq-BoundF} &
\min\limits_{y \in \E} \Bigl\{
F(y) + \frac{H + L_p}{(p + 1)!} \|y - x\|^{p + 1}
\Bigr\}.
\ea
$$
Therefore,
\beq \label{eq-GlobF}
\ba{rcl}
F(T(x)) - F^{*} & \leq &
\frac{H + L_p}{(p + 1)!}\|x - x^{*}\|^{p + 1}, \quad x \in \dom h,
\ea
\eeq
with $x^{*} \Def \arg\min\limits_{y \in \E} F(y)$, which
exists by our assumption. Denote by $D$ the maximal radius
of the initial level set of the objective, which we assume
to be finite:
$$
\ba{rcl}
D \;\; \Def \; \sup\limits_{x \in \dom h} \Bigl\{ \|x - x^{*}\|
:\; F(x) \leq F(x_0) \Bigr\}
\; & < & \; +\infty.
\ea
$$
Then, by monotonicity of the method~(\ref{met-RCTM}) and
by convexity we conclude
\beq\label{eq-ResG}
{1 \over D}\Bigl( F(x_{k + 1}) - F^* \Bigr)
\; \leq \; {1 \over D}\la F'(x_{k + 1}),
x_{k + 1} - x^{*} \ra  \;  \leq  \; \|F'(x_{k + 1})\|_{*}.
\eeq

In the general convex case, we can prove the global
sublinear rate of convergence of the Tensor Method of the
order $O({1 / k^p})$~\cite{nesterov2019implementable}. For
completeness of presentation, let us prove an extension of
this result onto the composite case.
\BT \label{th-SublR}
For the method~(\ref{met-RCTM}) with $H = pL_p$ we have
\beq\label{eq-SublR}
\ba{rcl}
F(x_{k}) - F^{*} & \leq & { (p + 1) (2p)^p \over p! }
\cdot {L_p D^{p + 1} \over (k - 1)^p}, \qquad k \geq 2.
\ea
\eeq
\ET
\proof Indeed, in view of~(\ref{eq-DecF}) and~(\ref{eq-ResG}),
we have for every $k \geq 0$
$$
\ba{rcl}
F(x_{k}) - F(x_{k + 1}) & \geq & \la F'(x_{k + 1}),
x_k - x_{k + 1} \ra\\
\\
& \refGE{eq-DecF} & \left( {p! \over (p+1)L_p}
\right)^{1 \over p} \cdot \| F'(x_{k + 1}) \|_*^{p+1 \over p}  \\
\\
& \refGE{eq-ResG} & \left( {p! \over (p+1)L_p D^{p + 1} }
\right)^{1 \over p}
\cdot \Bigl( F(x_{k + 1}) - F^* \Bigr)^{ p + 1 \over p }.
\ea
$$
Denoting $\delta_k = F(x_k) - F^*$ and
$C = \left( {p! \over (p+1) L_p D^{p + 1} }\right)^{1 \over p}$,
we obtain the following recurrence:
\beq\label{eq-Recurr}
\ba{rcl}
\delta_{k} - \delta_{k + 1} & \geq &
C \delta_{k + 1}^{p + 1 \over p}, \qquad k \geq 0,
\ea
\eeq
or for $\mu_k = C^p \delta_k \refLE{eq-GlobF} 1$, as follows:
$$
\ba{rcl}
\mu_{k} - \mu_{k + 1} & \geq & \mu_{k + 1}^{p + 1 \over p},
\qquad k \geq 0.
\ea
$$
Then, Lemma~1.1 from~\cite{grapiglia2017regularized}
provides us with the following guarantee:
$$
\ba{rcl}
\mu_{k} & \leq &
\Bigl(
\frac{p(1 + \mu_1^{1 / p})}{k - 1}
\Bigr)^p
\; \leq \;
\Bigl( \frac{2p}{k - 1} \Bigr)^p, \quad k \geq 2.
\ea
$$
Therefore,
$$
\ba{rcl}
\delta_k & = & {\mu_{k} \over C^p} \; \leq \;
\left({ 2p \over C (k - 1) }\right)^p \; = \;
{ (p + 1) (2p)^p \over p! } \cdot {L_p D^{p + 1} \over (k - 1)^p},
\qquad k \geq 2.
\ea
$$
\qed

For a given degree $q \geq 2$ of uniform convexity with
$\sigma_q > 0$, and for RCTM of order $p \geq q - 1$, let
us denote by~$\omega_{p, q}$ the following
\textit{condition number}:
$$
\ba{rcl}
\omega_{p, q} & \Def &
\frac{p + 1}{p!} \cdot
\Bigl(  \frac{q - 1}{q} \Bigr)^{q - 1}
\cdot \frac{L_p D^{p - q + 1}}{\sigma_q}.
\ea
$$

\BC\label{cor-Subl}
In order to achieve the region $\Q$ it is enough to perform
\beq \label{eq-Total1}
\Biggl\lceil
2p \cdot
\biggl(
{ q^{q} \over (q - 1)^{q - 1} }
\cdot
\omega_{p, q}^{\frac{p + 1}{p}}
\biggr)^{1 \over p - q + 1}
\Biggr\rceil + 2
\eeq
iterations of the method.
\EC
\proof
Plugging~(\ref{eq-RegF}) into~(\ref{eq-SublR}).
\qed

We can improve this estimate, knowing that the objective
is globally uniformly convex~(\ref{eq-UniC}). Then the
linear rate of convergence arises at the first state, till
the entering in the region~$\Q$.

\BT Let $\sigma_q > 0$ with $q \leq p + 1$.
Then for the method~(\ref{met-RCTM}) with $H = pL_p$, we
have
\beq\label{eq-LinR}
\ba{rcl}
F(x_{k}) - F^{*} & \leq & \exp\left( -{k \over 1 +
	\omega^{1/p}_{p, q}}
\right)
\cdot \bigl( F(x_{0}) - F^*\bigr), \qquad k \geq 1.
\ea
\eeq
Therefore, for a given $\varepsilon > 0$ to achieve
$F(x_K) - F^{*} \leq \varepsilon$, it is enough to set
\beq \label{eq-Total2}
\ba{rcl}
K & = & \left\lceil (1+\omega^{1/p}_{p,q}) \cdot
\log{\frac{F(x_0) - F^{*}}{\varepsilon}} \right\rceil + 1.
\ea
\eeq
\ET
\proof
Indeed, for every $k \geq 0$
$$
\ba{rcl}
F(x_{k}) - F(x_{k + 1}) & \geq & \la F'(x_{k + 1}),
x_k - x_{k + 1} \ra\\
\\
& \refGE{eq-DecF} & \left( {p! \over (p+1)L_p}
\right)^{1 \over p} \cdot \| F'(x_{k + 1}) \|_*^{p+1 \over p}  \\
\\
& = & \left( {p! \over (p+1)L_p} \right)^{1 \over p}
\cdot \| F'(x_{k + 1}) \|_*^{p - q + 1 \over p}
\cdot \| F'(x_{k + 1}) \|_*^{q \over p} \\
\\
& \stackrel{(\ref{eq-ResG}),(\ref{eq-DecFG})}{\geq} &
\left( {p! \over p + 1} \cdot
{ \sigma_q \over L_p D^{p - q + 1}} \right)^{1 \over p}
\cdot \left({ q \over q - 1 }\right)^{q - 1 \over p}
\cdot \Bigl( F(x_{k + 1}) - F^* \Bigr) \\
\\
& = &
\left( \frac{1}{\omega_{p, q}}  \right)^{1 \over p}
\cdot \Bigl( F(x_{k + 1}) - F^* \Bigr).
\ea
$$
Denoting $\delta_k = F(x_k) - F^{*}$, we obtain
$$
\ba{rcl}
\delta_{k + 1} & \leq & {\omega^{1/p}_{p,q} \over 1 +
	\omega^{1/p}_{p,q}} \cdot \delta_k \; \leq \; \exp \left(
- {1 \over 1 + \omega^{1/p}_{p,q}} \right) \cdot \delta_k,
\qquad k \geq 1.
\ea
$$
\qed

We see that, for RCTM with $p \geq 2$ minimizing the
uniformly convex objective of degree $q \leq p + 1$, the
condition number $\omega^{1/p}_{p, q}$ is the main factor
in the global complexity estimates~(\ref{eq-Total1})
and~(\ref{eq-Total2}). Since in general this number may be
arbitrarily big, complexity estimate $\tilde{O}(\omega_{p,
	q}^{1 / p})$ in (\ref{eq-Total2}) is much better than the
estimate $O(\omega_{p, q}^{(p + 1) / (p(p - q + 1))})$ in
(\ref{eq-Total1}) because of relation ${ p + 1 \over p - q
	+ 1} \geq 1$.

These global bounds can be improved, by using the
\textit{universal}~\cite{doikov2019minimizing,grapiglia2019tensor}
and the
\textit{accelerated}~\cite{nesterov2008accelerating,grapiglia2019accelerated,grapiglia2019tensor,gasnikov2019optimal,song2019towards}
high-order schemes.

High-order tensor methods for minimizing
the gradient norm were developed in
\cite{dvurechensky2019near}.
These methods achieve near-optimal global convergence rates,
and can be used for coming into the region~$\G$~\eqref{eq-RegG}.
Note, that for the composite minimization problems,
some modification of these methods is required,
which ensures minimization of the \textit{subgradient} norm.

Finally, let us mention some recent results~\cite{nesterov2020superfast,kamzolov2020near},
where it was shown that
a proper implementation of the third-order schemes by
second-order oracle may lead to a significant acceleration of the methods.
However, the relation of these techniques to the local convergence
needs further investigations.

\section{Application to proximal methods}
\label{sc-Prox}
\SetEQ

Let us discuss now a general approach, which uses the
local convergence of the methods for justifying the global
performance of proximal iterations.

The proximal method~\cite{rockafellar1976monotone} is
one of the classical methods in
theoretical optimization.
Every step of the method for solving problem~(\ref{prob-Main})
is a minimization of the regularized objective:
\beq \label{Prox-Subprob}
\ba{rcl}
x_{k + 1} & = & \arg\min\limits_{x \in \E}
\Bigl\{
a_{k + 1} F(x) + \frac{1}{2}\|x - x_k\|^2
\Bigr\}, \qquad k \geq 0,
\ea
\eeq
where $\{ a_k \}_{k \geq 1}$ is a sequence of positive
coefficients, related to the iteration counter.

Of course, in general, we can hope only to solve
subproblem~(\ref{Prox-Subprob}) inexactly. The questions
of practical implementations and possible
generalizations of the proximal method, are still in the
area of intensive research (see, for example
\cite{guler1991convergence,solodov2001unified,schmidt2011convergence,salzo2012inexact}).

One simple observation on the
subproblem~(\ref{Prox-Subprob}) is that it is $1$-strongly
convex. Therefore, if we would be able to pick an initial
point from the region of superlinear
convergence~(\ref{eq-RegF}) or~(\ref{eq-RegG}), we could
minimize it very quickly by RCTM of degree $p \geq 2$ up
to arbitrary accuracy. In this section, we are going to
investigate this approach. For the resulting scheme, we
will prove the global rate of convergence of the order
$\tilde{O}(1 / k^{p + 1 \over 2})$.

Denote by $\Phi_{k + 1}$ the regularized objective
from~(\ref{Prox-Subprob}):
$$
\ba{rcl}
\Phi_{k + 1}(x) & \Def &
a_{k + 1} F(x) + \frac{1}{2}\|x - x_k\|^2
\; = \;
a_{k + 1} f(x) + \frac{1}{2}\|x - x_k\|^2
+ a_{k + 1} h(x).
\ea
$$
We fix a sequences of accuracies $\{\delta_k\}_{k \geq 1}$
and relax the assumption on exact minimization
in~(\ref{Prox-Subprob}). Now, at every step we need to
find a point $x_{k + 1}$ and corresponding subgradient
vector $g_{k + 1} \in \partial \Phi_{k + 1}(x_{k + 1})$
with bounded norm:
\beq \label{RelaxedMin}
\ba{rcl}
\|g_{k + 1}\|_{*} & \leq & \delta_{k + 1}.
\ea
\eeq
Denote
$$
\ba{rcl}
F'(x_{k + 1}) & \Def &
\frac{1}{a_{k + 1}}( g_{k + 1} - B(x_{k + 1} - x_k))
\; \in \; \partial F(x_{k + 1}).
\ea
$$

The following global convergence result holds
for the general proximal method with inexact
minimization criterion~(\ref{RelaxedMin}).

\BT \label{th-InexProx}
Assume that there exist a minimum $x^{*} \in \dom h$ of
the problem~(\ref{prob-Main}). Then, for any $k \geq 1$,
we have
\beq \label{eq-InexProx}
\ba{rcl}
\sum\limits_{i = 1}^k a_i(F(x_i) - F^{*})
+ \frac{1}{2}\sum\limits_{i = 1}^k a_i^2 \|F'(x_i)\|_{*}^2
+ \frac{1}{2}\|x_k - x^{*}\|^2
& \leq & R_k(\delta),
\ea
\eeq
where
$$
\ba{rcl}
R_k(\delta) & \Def &
\frac{1}{2}\left(
\|x_0 - x^{*}\| + \sum\limits_{i = 1}^k \delta_i
\right)^2.
\ea
$$
\ET
\proof
First, let us prove that for all $k \geq 0$ and for every
$x \in \dom h$, we have
\beq \label{InductionCondition}
\ba{rcl}
\frac{1}{2}\|x_0 - x\|^2
+ \sum\limits_{i = 1}^k a_i F(x) 
& \geq & 
\frac{1}{2}\|x_k - x\|^2 + C_k(x),
\ea
\eeq
where
$$
\ba{rcl}
C_k(x) & \Def &
\sum\limits_{i = 1}^k \left( a_i F(x_i) +
\frac{a_i^2}{2} \|F'(x_i)\|_{*}^2 + \la g_i, x - x_{i - 1}
\ra - \frac{\delta_i^2}{2} \right).
\ea
$$
This is obviously true for $k = 0$. Let it hold for some
$k \geq 0$. Consider the step number $k + 1$ of the
inexact proximal method.

By condition~(\ref{RelaxedMin}), we have
$$
\ba{rcl}
\| a_{k + 1} F'(x_{k + 1}) + B(x_{k + 1} - x_k) \|_{*}^2 &
\leq & \delta_{k + 1}^2.
\ea
$$
Equivalently,
\beq \label{ProxOneStep}
\ba{cl}
& \la a_{k + 1} F'(x_{k + 1}), x_k - x_{k + 1} \ra \\
\\
& \; \geq \;
\frac{a_{k + 1}^2}{2}\|F'(x_{k + 1})\|_{*}^2
+ \frac{1}{2}\|x_{k + 1} - x_k\|^2
- \frac{\delta_{k + 1}^2}{2}.
\ea
\eeq
Therefore, using the inductive assumption and strong
convexity of $\Phi_{k + 1}(\cdot)$, we conclude
$$
\ba{rl}
& \frac{1}{2}\|x_0 - x\|^2 + \sum\limits_{i = 1}^{k + 1}
a_i F(x) \; = \; \frac{1}{2}\|x_0 - x\|^2 + \sum\limits_{i
	= 1}^k a_i F(x) + a_{k + 1} F(x)
\\
\\
& \; \refGE{InductionCondition} \; 
\Phi_{k + 1}(x) + C_k(x) \\
\\
& \;\;\, \geq \;\;\, 
\Phi_{k + 1}(x_{k + 1}) + \la g_{k + 1}, x - x_{k + 1} \ra
+ \frac{1}{2}\|x_{k + 1} - x\|^2 + C_k(x) \\
\\
& \;\;\, = \;\;\,  a_{k + 1} F(x_{k + 1}) + \frac{1}{2}\|x_{k + 1} -
x_k\|^2
+ \la g_{k + 1}, x_k - x_{k + 1} \ra   \\
\\
& \;\;\qquad + \;\;\, \la g_{k + 1}, x - x_k \ra
+ \frac{1}{2}\|x_{k + 1} - x\|^2 + C_k(x) \\
\\
& \;\;\, = \;\;\,  
a_{k + 1} F(x_{k + 1}) + \la a_{k + 1} F'(x_{k + 1}),
x_k - x_{k + 1} \ra
- \frac{1}{2}\|x_{k + 1} - x_k\|^2  \\
\\
& \;\;\qquad + \;\;\, \la g_{k + 1}, x - x_k \ra
+ \frac{1}{2}\|x_{k + 1} - x\|^2 + C_k(x) \\
\\
& \; \refGE{ProxOneStep} \; a_{k + 1} F(x_{k + 1}) + \frac{a_{k
		+ 1}^2}{2}\|F'(x_{k + 1})\|_{*}^2
- \frac{\delta_{k + 1}^2}{2} \\
\\
& \;\;\qquad + \;\;\, 
\la g_{k + 1}, x - x_k \ra + \frac{1}{2}\|x_{k + 1} - x\|^2 + C_k(x) \\
\\
& \;\;\, = \;\;\,
\frac{1}{2}\|x_{k + 1} - x\|^2 + C_{k + 1}(x).
\ea
$$
Thus, inequality~(\ref{InductionCondition}) is valid for
all $k \geq 0$.

Now, plugging $x \equiv x^{*}$ into~(\ref{InductionCondition}),
we have
\beq \label{RecurProx}
\ba{cl}
&  \sum\limits_{i = 1}^k a_i (F(x_i) - F^{*})
+ \frac{1}{2}\sum\limits_{i = 1}^k a_i^2 \|F'(x_i)\|_{*}^2
+ \frac{1}{2}\|x_k - x^{*}\|^2 \\
\\
& \;\; \, \leq \;\;\,
\frac{1}{2}\|x_0 - x^{*}\|^2
+ \frac{1}{2}\sum\limits_{i = 1}^k \delta_i^2
+ \sum\limits_{i = 1}^k \la g_i, x_{i - 1} - x^{*} \ra \\
\\
& \; \refLE{RelaxedMin} \;
\frac{1}{2}\|x_0 - x^{*}\|^2
+ \frac{1}{2}\sum\limits_{i = 1}^k \delta_i^2
+ \sum\limits_{i = 1}^k \delta_i \|x_{i - 1} - x^{*} \|
\quad \Def \quad \alpha_k.
\ea
\eeq
In order to finish the proof, it is enough to show that
$\alpha_k \leq R_k(\delta)$.

Indeed,
$$
\ba{rcl}
\alpha_{k + 1} & = &
\alpha_k + \frac{1}{2} \delta_{k + 1}^2
+ \delta_{k + 1} \|x_k - x^{*}\| \\
\\
& \refLE{RecurProx} &
\alpha_k + \frac{1}{2}\delta_{k + 1}^2
+ \delta_{k + 1} \sqrt{2 \alpha_k} \\
\\
& = &
\left( \sqrt{\alpha_k}
+ \frac{1}{\sqrt{2}}\delta_{k + 1} \right)^2.
\ea
$$
Therefore,
$$
\ba{rcl}
\sqrt{\alpha_k} & \leq &
\sqrt{\alpha_{k - 1}} + \frac{1}{\sqrt{2}}\delta_{k}
\; \leq \; \dots \; \leq \;
\sqrt{\alpha_0} + \frac{1}{\sqrt{2}}\sum\limits_{i = 1}^k \delta_i \\
\\
& = &
\frac{1}{\sqrt{2}}\left( \|x_0 - x^{*}\|
+ \sum\limits_{i = 1}^k \delta_i \right)
\; = \; \sqrt{R_k(\delta)}.
\ea
$$
\qed

Now, we are ready to use the result on the local
superlinear convergence of RCTM in the norm of subgradient
(Theorem~\ref{th-RateG}), in order to minimize $\Phi_{k +
	1}(\cdot)$ at every step of inexact proximal method.

Note that
$$
\ba{rcl}
\partial \Phi_{k + 1}(x) & = & a_{k + 1} \partial F(x) + B(x - x_k),
\ea
$$
and it is natural to start minimization process from the
previous point $x_k$, for which $\partial \Phi_{k +
	1}(x_k) = a_{k + 1} \partial F(x_k)$. Let us also notice,
that the Lipschitz constant of the $p$th derivative ($p
\geq 2$) of the smooth part of~$\Phi_{k + 1}$ is $a_{k +
	1} L_p$.

Using our previous notation, one step of RCTM can be
written as follows:
$$
\ba{cl}
&T_H(\Phi_{k + 1}, z) \\
\\
& \;\;\, \Def \;\;\,
\arg\min\limits_{y \in \E}
\Bigl\{
a_{k + 1} \Omega_{p}(f, z; y) + \frac{H}{(p + 1)!}\|y - z\|^{p + 1}
+ a_{k + 1}h(y)
+ \frac{1}{2}\|y - x_k\|^2
\Bigr\},
\ea
$$
where $H = a_{k + 1}pL_p$. Then, a sufficient condition
for $z = x_k$ to be in the region of superlinear
convergence~\eqref{eq-RegG} is
$$
\ba{rcl}
a_{k + 1} \| F'(x_k) \|_*
& \leq &
\left(
p! \over a_{k + 1} (p + 1)  L_p
\right)^{1 \over p - 1},
\ea
$$
or, equivalently
$$
\ba{rcl}
a_{k + 1} & \leq &
\left({1 \over \|F'(x_k)\|_{*} }\right)^{p - 1 \over p}
\left({ p! \over (p + 1) L_p  }\right)^{1 \over p}.
\ea
$$
To be sure that $x_k$ is strictly inside the region, we
can pick:
\beq \label{ak-Choice}
\boxed{
	\ba{rcl}
	a_{k + 1} & = &
	\left({1 \over 2 \| F'(x_k)\|_{*}} \right)^{p - 1 \over p}
	\left(
	p! \over (p + 1) L_p
	\right)^{1 \over p}
	\ea
}
\eeq
Note, that this rule requires fixing an initial
subgradient $F'(x_0) \in \partial F(x_0)$, in order to
choose $a_1$.

Finally, we apply the following steps:
\beq\label{met-RCTM-2}
\ba{rcl}
z_0 \; = \; x_k, \quad z_{t+1} & = & T_{H}(\Phi_{k + 1}, z_t), \quad
t \geq 0.
\ea
\eeq
We can estimate the required number of these iterations as follows.
\BL
At every iteration $k \geq 0$ of the inexact proximal
method, in order to achieve $\| \Phi'_{k + 1}(z_t) \|_{*}
\leq \delta_{k + 1}$, it is enough to perform
\beq \label{eq-LogLog}
\ba{rcl}
t_k &=&
\biggl\lceil \frac{1}{\log_2 p} \cdot  \log_2 \log_2
\left(
\frac{2 D_k(\delta) }{ \delta_{k + 1}}
\right) \biggr\rceil
\ea
\eeq
steps of RCTM~\eqref{met-RCTM-2},
where
$$
\ba{rcl}
D_k(\delta) & \Def &
\max \biggl\{ 
\|x_0 - x^{*}\| + \sum\limits_{i = 1}^k \delta_i, 
\Bigl(  \frac{p! \|F'(x_0)\|_{*} }{(p + 1)L_p2^{p - 1}}  \Bigr)^{1 \over p}
\biggr\}
\ea
$$
\EL
\proof
According to~\eqref{eq-RateG}, one step of RCTM~\eqref{met-RCTM-2}
provides us with the following guarantee 
in terms of the subgradients of our objective $\Phi_{k + 1}(\cdot)$:
\beq \label{step-RCTM-2}
\ba{rcl}
\| \Phi'_{k + 1}(z_t) \|_{*}
& \leq & 
\frac{a_{k + 1} (p + 1) L_p}{p!} \| \Phi'_{k + 1}(z_{t - 1}) \|_{*}^p,
\ea
\eeq
where we used in~\eqref{eq-RateG} the values $q = 2$, $\sigma_q = 1$,
$a_{k + 1} L_p$ for the Lipschitz constant of the $p$th derivative of the smooth part of $\Phi_{k + 1}$,
and $H = a_{k + 1}pL_p$.

Denote
$\beta \equiv \left( { a_{k + 1}(p + 1)L_p \over p!  } \right)^{1 \over p - 1}
\refEQ{ak-Choice}
\left( {  (p + 1) L_p \over 2 \cdot p! \cdot \|F'(x_k)\|_* } \right)^{1 \over p}$.
Then, from~\eqref{step-RCTM-2} we have
\beq \label{Sublin-Conv}
\ba{rcl}
\beta \| \Phi'_{k + 1}(z_t) \|_{*}  & \leq &
\bigl(\beta \| \Phi'_{k + 1}(z_{t - 1}) \|_{*}\bigr)^{p} \\
\\
& \leq & \dots \;\; \leq \;\;
\bigl(\beta \| \Phi'_{k + 1}(z_0) \|_{*}\bigr)^{p^t} \\
\\
& = &
(\beta a_{k + 1}\|F'(x_k)\|_{*})^{p^t} \\
\\
& = &
\left(
a_{k + 1}^{p \over p - 1}
\left({ (p + 1) L_p \over p!  }\right)^{1 \over p - 1}
\|F'(x_k)\|_{*}
\right)^{p^t} \\
\\
& \refEQ{ak-Choice} & \left({1 \over 2}\right)^{p^t}.
\ea
\eeq
Therefore, for
\beq \label{eq-LogLog-2}
\ba{rcl}
t & \geq &
\log_p \log_2 \left( \frac{1}{\beta \delta_{k + 1}} \right)
\; = \;
\frac{1}{\log_2 p} \cdot \log_2 \log_2
\left(
\frac{1}{ \delta_{k + 1}}
\left( {
	2 \cdot p! \cdot \| F'(x_k) \|_* \over (p + 1) L_p
} \right)^{1 \over p} \right),
\ea
\eeq
it holds $\| \Phi'_{k + 1}(z_t)\|_{*} \leq \delta_{k +
	1}$. To finish the proof, let us estimate $\| F'(x_k)
\|_{*}$ from above. We have
\beq \label{eq-GBound}
\ba{rcl}
2^{3p - 2 \over p} \left(  \frac{(p + 1)L_p}{p!}  \right)^{2 \over p} R_k(\delta)
& \refGE{eq-InexProx} &
2^{2(p - 1) \over p} \left(  \frac{(p + 1)L_p}{p!}  \right)^{2 \over p}
\sum\limits_{i = 1}^k a_i^2 \|F'(x_i)\|_{*}^2 \\
\\
& \refEQ{ak-Choice} &
\sum\limits_{i = 1}^k \|F'(x_{i - 1})\|_{*}^{2(1 - p) \over p} \|F'(x_i)\|_{*}^2.
\ea
\eeq
Thus, for every $1 \leq i \leq k$ it holds
\beq \label{eq-GBound2}
\ba{rcl}
\|F'(x_i)\|_{*} & \refLE{eq-GBound} & \| F'(x_{i - 1})\|_{*}^{\rho}
\cdot \mathcal{D},
\ea
\eeq
with
$\mathcal{D} \equiv R_k^{1/2}(\delta)
\left( \frac{(p + 1) L_p}{p!}  \right)^{1 \over p} 2^{3p - 2 \over 2p}$,
and $\rho \equiv \frac{p - 1}{p}$.
Therefore,
$$
\ba{rcl}
\|F'(x_k)\|_{*} & \refLE{eq-GBound2} & \|F'(x_0)\|_{*}^{\rho^k}
\cdot \mathcal{D}^{1 + \rho + \rho^2 + \dots + \rho^{k - 1}} \\
\\
& = & \|F'(x_0)\|_{*} \cdot
\Bigl(  \|F'(x_0)\|_{*}^{\rho^k - 1} \cdot \mathcal{D}^{\frac{\; \;1 - \rho^k}{1 - \rho}}  \Bigr) \\
\\
& = & \|F'(x_0)\|_{*} \cdot \left(
\frac{\mathcal{D}^{p}}{\|F'(x_0)\|_{*}}  \right)^{1 - \rho^k}
\; \leq \;
\| F'(x_0) \|_{*} \cdot \max \bigl\{  \frac{\mathcal{D}^p}{\|F'(x_0)\|_{*}}, 1  \bigr\} \\
\\
& = & \max \biggl\{
\frac{(p + 1) L_p 2^{p - 1}}{p!}
\Bigl( \|x_0 - x^{*}\| + \sum\limits_{i = 1}^k \delta_i
\Bigr)^p, \; \|F'(x_0)\|_{*} \biggr\}.
\ea
$$
Substitution of this bound into~\eqref{eq-LogLog-2}
gives~\eqref{eq-LogLog}.
\qed

Let us prove now the rate of convergence for the outer
iterations. This is a direct consequence of
Theorem~\ref{th-InexProx} and the choice~\eqref{ak-Choice}
of the coefficients $\{ a_{k} \}_{k \geq 1}$.

\BL Let for a given $\varepsilon > 0$,
\beq \label{eq-epsLBound}
\ba{rcl}
F(x_k) - F^{*} & \geq & \varepsilon, \qquad 1 \leq k \leq K.
\ea
\eeq
Then for every $1 \leq k \leq K$, we have
\beq \label{eq-cProx}
\ba{rcl}
F(\bar{x}_k) - F^{*} & \leq &
\frac{L_p \left(
	\|x_0 - x^{*} \| + \sum_{i = 1}^k \delta_i
	\right)^{p + 1}
}{k^{p + 1 \over 2}}
\frac{(p + 1) 2^{p - 2} V_k(\varepsilon) }{ p!},
\ea
\eeq
where
$\bar{x}_k \Def \frac{\sum_{i = 1}^k a_i x_i}{\sum_{i = 1}^k a_i}$, and
$V_k(\varepsilon) \Def \left( \frac{\|F'(x_0)\|_{*} \cdot ( \|x_0 - x^{*}\|
	+ \sum_{i = 1}^k \delta_i  )}{\varepsilon}
\right)^{p - 1 \over k}$.
\EL
\proof
Using the inequality between the arithmetic and geometric
means, we obtain
\beq \label{ak-rate}
\ba{rcl}
R_{k}(\delta)
& \refGE{eq-InexProx} &
\frac{1}{2}\sum\limits_{i = 1}^k a_i^2 \|F'(x_i)\|_*^2
\; \refEQ{ak-Choice} \;
\frac{1}{8}
\left( \frac{p!}{(p + 1)L_p}
\right)^{2 \over p - 1}
\sum\limits_{i = 1}^k
\frac{a_i^2}{a_{i + 1}^{2p \over p - 1}} \\
\\
& \geq &
\frac{k}{8}
\left( \frac{p!}{(p + 1)L_p}
\right)^{2 \over p - 1}
\left(
\prod\limits_{i = 1}^k
\frac{a_i^2}{a_{i + 1}^{2p \over p - 1}}
\right)^{1 \over k} \\
\\
& = &
\frac{k}{8}
\left( \frac{p!}{(p + 1)L_p}
\right)^{2 \over p - 1}
\left(
\frac{a_1}{a_{k + 1}}
\right)^{2p \over (p - 1)k}
\left( \prod\limits_{i = 1}^k a_i \right)^{-2 \over (p - 1)k} \\
\\
& \geq &
\frac{k^{p + 1 \over p - 1}}{8}
\left( \frac{p!}{(p + 1)L_p}
\right)^{2 \over p - 1}
\left(
\frac{a_1}{a_{k + 1}}
\right)^{2p \over (p - 1)k}
\left( \sum\limits_{i = 1}^k a_i \right)^{-2 \over p - 1}.
\ea
\eeq
Therefore,
$$
\ba{rcl}
F(\bar{x}_k) - F^{*}
& \leq &
\frac{1}{\sum\limits_{i = 1}^k a_i}
\sum\limits_{i = 1}^k a_i (F(x_i) - F^{*})
\; \refLE{eq-InexProx} \;
\frac{R_k(\delta)}{\sum\limits_{i = 1}^k a_i} \\
\\
& \refLE{ak-rate} &
\frac{ R_k(\delta)^{p + 1 \over 2}  }{k^{p + 1 \over 2}}
\frac{(p + 1) L_p}{p!}
\left( \frac{a_{k + 1}}{a_1} \right)^{p \over k}
8^{p - 1 \over 2} \\
\\
& = &
\frac{L_p \left(
	\|x_0 - x^{*} \| + \sum_{i = 1}^k \delta_i
	\right)^{p + 1}
}{k^{p + 1 \over 2}}
\frac{(p + 1) 2^{p - 2} }{ p!}
\left( \frac{\|F'(x_0)\|_{*}}{\|F'(x_k)\|_{*}} \right)^{p - 1 \over k},
\ea
$$
where the first inequality holds by convexity.
At the same time, we have
$$
\ba{rcl}
\|F'(x_k)\|_{*} & \geq & \frac{\la F'(x_k), x_k - x^{*} \ra}{\|x_k - x^{*}\|}
\; \geq \; \frac{F(x_k) - F^{*}}{\|x_k - x^{*}\|} \\
\\
& \refGE{eq-epsLBound} & \frac{\varepsilon}{\|x_k - x^{*}\|}
\; \refGE{eq-InexProx} \; \frac{\varepsilon}{\|x_0 - x^{*}\| + \sum_{i = 1}^k \delta_i }.
\ea
$$
Thus, $\left( \frac{\|F'(x_0)\|_{*}}{\|F'(x_k)\|_{*}} \right)^{p - 1 \over k} \leq V_k(\varepsilon)$
and we obtain~\eqref{eq-cProx}.
\qed

\BR
Note that
$\bigl(\frac{1}{\varepsilon}\bigr)^{p - 1 \over k}
= \exp\bigl( {p - 1 \over k} \ln {1 \over \varepsilon} \bigr)$.
Therefore after
$k = O\left( \ln {1 \over \varepsilon}\right)$ iterations, the factor $V_k(\varepsilon)$
is bounded by an absolute constant.
\ER

Since the local convergence of RCTM is very
fast~\eqref{eq-LogLog}, we can choose the inner
accuracies~$\{ \delta_i \}_{i \geq 1}$ small enough, to
have the right hand side of~\eqref{eq-cProx} being of the
order $\tilde{O}(1 / k^{p + 1 \over 2})$. Let us present a
precise statement.

\BT
Let $\delta_k \equiv \frac{c}{k^s}$ for fixed absolute
constants $c > 0$ and $s > 1$. Let for a given
$\varepsilon > 0$, we have
$$
\ba{rcl}
F(x_k) - F^{*} & \geq & \varepsilon, \qquad
1 \leq k \leq K.
\ea
$$
Then, for every $k$ such that $\ln \frac{\|F'(x_0)\|_{*}
	R}{ \varepsilon}  \leq k \leq K$, we get
\beq \label{eq-PrConv}
\ba{rcl}
F(\bar{x}_k) - F^{*} & \leq &
\frac{L_p R^{p + 1}}{k^{p + 1 \over 2}} \frac{(p + 1) 2^{p - 2} \exp(p - 1)}{p!},
\ea
\eeq
where
$$
\ba{rcl}
R & \Def & \|x_0 - x^{*}\| + \frac{cs}{s - 1}.
\ea
$$
The total number of oracle calls $N_k$ during the first
$k$ iterations is bounded as follows:
$$
\ba{rcl}
N_k & \leq & k \cdot \Bigl( 1 + \frac{1}{\log_2 p}  \log_2 \log_2 \frac{2D k^s }{c} \Bigr),
\ea
$$
where
$$
\ba{rcl}
D & \Def & \max \biggr\{
R, \,
\Bigl(  \frac{p! \|F'(x_0)\|_{*} }{(p + 1)L_p2^{p - 1}}  \Bigr)^{1 \over p}
\biggl\}.
\ea
$$
\ET
\proof
Indeed,
$$
\ba{rcl}
\sum\limits_{i = 1}^k \delta_i & = &
c\Bigl(1 + \sum\limits_{i = 2}^k \frac{1}{i^s} \Bigr)
\; \; \leq \; \;
c\Bigl(1 + \int\limits_1^k \frac{dx}{x^{s}} \Bigr)
\; \; = \; \;
c\Bigl(1 - \frac{1}{s - 1} \int\limits_1^k dx^{-(s - 1)} \Bigr) \\
\\
& = &
c\Bigl(1  - \frac{k^{-(s - 1)}}{s - 1} + \frac{1}{s - 1} \Bigr)
\; \; \leq \; \; \frac{cs}{s - 1}.
\ea
$$
Thus, we obtain~\eqref{eq-PrConv} directly from
the bound~\eqref{eq-cProx}, and by the fact that
$$
\ba{rcl}
V_k(\varepsilon) & \equiv &
\Bigl(  \frac{\| F'(x_0) \|_{*} R}{\varepsilon}  \Bigr)^{\frac{p - 1}{k}}
\; = \;
\exp\Bigl( \frac{p - 1}{k} \log \frac{\| F'(x_0) \|_{*} R}{\varepsilon}  \Bigr) \\
\\
& \leq & \exp(p - 1),
\ea
$$
when $k \geq \ln \frac{\| F'(x_0) \|_{*} R }{ \varepsilon} $.

Finally,
$$
\ba{rcl}
N_k & \refLE{eq-LogLog} &
\sum\limits_{i = 1}^k \left\lceil  \frac{1}{\log_2 p}
\log_2 \log_2 \frac{2 D }{\delta_i} \right\rceil
\; \leq \;
k +  \frac{1}{\log_2 p} \sum\limits_{i = 1}^k \log_2 \log_2 \frac{2Di^s}{c} \\
\\
& \leq & k + \frac{1}{\log_2 p} \sum\limits_{i = 1}^k \log_2 \log_2 \frac{2Dk^s}{c}
\; = \;
k \cdot \Bigl(1 + \frac{1}{\log_2 p} \log_2 \log_2 \frac{2Dk^s}{c} \Bigr).
\ea
$$
\qed

Note that we were able to justify the global performance
of the scheme, using only the local convergence results
for the inner method. It is interesting to compare our
approach with the recent results on the path-following
second-order methods \cite{dvurechensky2018global}.

We can drop the logarithmic components in the complexity
bounds by using the \textit{hybrid proximal methods}
(see~\cite{monteiro2010complexity}
and~\cite{marques2019iteration}), where at each iteration
only one step of RCTM is performed. The resulting rate of
convergence there is $O(1 / k^{p + 1 \over 2})$, without
any extra logarithmic factors. However, this rate is worse
than the rate $O(1 / k^p)$ provided by the
Theorem~\ref{th-SublR} for the primal iterations of
RCTM~\eqref{met-RCTM}.

\section*{Acknowledgements}

We are very thankful to anonymous referees
for valuable comments that improved the initial version of this paper.

%\bibliographystyle{spmpsci.bst}
%\bibliography{bibliography}

\end{document}